\magnification=1200

\noindent {\bf A stabilization theorem for Hermitian forms and applications to holomorphic mappings}

David W. Catlin and John P. D'Angelo

Catlin

Dept. of Mathematics

Purdue University

W. Lafayette IN 47907

D'Angelo

Dept. of Mathematics

University of Illinois

Urbana IL 61801

\noindent {\bf Introduction}.

We consider positivity conditions both for real-valued functions of several complex variables 
and for Hermitian forms. We prove a stabilization theorem relating these two notions,
and give some applications to proper mappings between balls in different dimensions. 
The technique of proof relies on the simple expression for the Bergman kernel function
for the unit ball and elementary facts about Hilbert spaces. Our main result generalizes to Hermitian forms 
a theorem
proved by Polya [HLP] for homogeneous real polynomials, which was obtained in conjunction with Hilbert's
seventeenth problem. See [H] and [R] for generalizations of Polya's theorem of a completely different kind.
The flavor of our applications is also completely different.

We begin by describing our main result. 
Let ${\bf {C^n}}$ denote complex 
Euclidean space 
of $n$ dimensions, with complex Euclidean squared norm $||z||^2 = \sum_{j=1}^n |z_j|^2$.
Suppose that $ f:{\bf C^n} \rightarrow {\bf R} $ is a polynomial in the variables
 $z$ and ${\overline z}$, and that
 it is homogeneous of the
same degree $m$ in each of these variables. We write 
$f(z,{\overline z}) = \sum_{|\alpha|=m} \sum_{|\beta|=m} c_{\alpha \beta} z^\alpha {\overline z}^\beta$.
 The
condition that $f$ take only real values is equivalent to 
$c_{\alpha \beta} = {\overline {c_{\beta \alpha}}}$. 
We call the Hermitian matrix $(c_{\alpha \beta})$ the underlying matrix of coefficients.
Then we have the following conclusion. The function $f$ achieves a positive minimum value on the
unit sphere if and only if there is an integer $d$ so that the matrix $(E_{\mu \nu})$ is
 positive definite, where $(E_{\mu \nu})$ is 
the underlying matrix of coefficients for the function $||z||^{2d} f(z,{\overline z})$, that is 
$$||z||^{2d} f(z,{\overline z}) = \sum E_{\mu \nu} z^\mu {\overline z}^\nu .$$
Consequently there is a homogeneous holomorphic vector-valued polynomial $g$ such that

$$||z||^{2d} f(z,{\overline z}) = ||g(z)||^2 $$

When $(c_{\alpha \beta})$ is itself positive definite, then
$f$ must be positive on the sphere; this is the case $d=0$.
In general the smallest possible value for the integer $d$ depends on the original 
underlying matrix of coefficients. We note that, once the form is positive for some $d$, it remains
 positive for all larger values, and this suggests the name ``stabilization''.
In case the original underlying matrix of coefficients is diagonal, this 
theorem implies the classical theorem of Polya in the real case. Even in the diagonal case the smallest possible 
integer $d$ can be arbitrarily large. 
See Example 2 and Remark 1. 

It is elementary that when the underlying matrix of coefficients is positive
definite the underlying polynomial is the squared norm of a holomorphic mapping. We conclude that,
although $f(z,\overline {z})$ may not be the squared norm of a holomorphic mapping,
 for sufficiently large integers $d$,
$||z||^{2d} f(z, \overline{z})$ is such a squared norm $||g(z)||^2$. Furthermore this function
vanishes only at the origin.

We next indicate how this result applies to holomorphic mappings.
In Theorem 3 we let  $f$ be an arbitrary 
polynomial which is positive on the sphere; using Theorem 1 we construct a holomorphic
polynomial mapping $g$ (finite-dimensionally valued) 
such that $f(z,{\overline z}) = ||g(z)||^2$ on the sphere.
Theorem 3 easily implies Theorem 2. In Theorem 2, we consider a holomorphic polynomial mapping $p$ 
whose Euclidean norm on the closed ball in ${\bf C^n}$ is everywhere less than unity. Then we can
find a holomorphic polynomial mapping
$g$ such that $p\oplus g$ maps the sphere to some sphere; in other words,
$p$ equals some of the components of a proper holomorphic polynomial mapping between balls. 
We must allow the target dimension to be arbitrarily large. This cannot be improved, because 
even in the homogeneous case, the (minimum) embedding dimension
of $p \oplus g$ depends on the integer $d$ from Theorem 1. In Theorem 4, we
 give a result on allowable denominators
for rational proper mappings between balls.

We also interpret Theorem 1 in terms of tensor products of mappings.
If $h$ and $g$ are vector valued functions, we define their tensor product $h\otimes g$ to be the mapping
whose components are all possible products $h_j g_k$ of the components of $h$ and $g$, in some determined order.
Then $||z||^{2d} =||H_{d}(z)||^2$ for a certain holomorphic mapping $H_d$, namely the $d$-fold 
tensor product of the identity mapping with
itself. See [D1] for many uses of this mapping.  We observe that we may always write a real-valued 
polynomial $p$ as 
$p(z,{\overline z}) = ||P(z)||^2 - ||N(z)||^2 $ for holomorphic mappings $P$ and $N$. 
Theorem 1 then implies that
$p$ is positive on the sphere if and only if there is an integer $d$ so that we can write 

$$ H_d \otimes N = L( H_d \otimes  P)$$ 
where $L$ is a linear transformation in the unit ball of the space of linear transformations.
 See equivalence 5 of Theorem 1
for the precise statement.

It is possible to phrase our result in terms of Hermitian metrics
on powers of the hyperplane section bundle of complex projective space, but we do
not do this here. Because of its interpretation involving tensor products,
our statement might be related to the famous theorem of Kodaira on the embedding
 of Hodge manifolds into complex
projective spaces, where one takes sufficiently high tensor powers of a positive line bundle. See [GH].
We do not pursue this here.
At the close of the paper we briefly
discuss the geometric meaning of the main theorem in terms of the Veronese mapping. 

We close the introduction by giving another interpretation.
Letting the $n$-torus act on the sphere, we see that the positivity of $f$ as a function shows that 
certain trigonometric polynomials are positive. By the easy direction of 
Bochner's theorem characterizing functions of positive type,
all the matrices in a certain family 
must be positive definite. We relate the matrix entries of the 
underlying form for $||z||^{2d} f(z, {\overline z})$ to the matrices in this family when $d$
is sufficiently large. This argument
relies on the elementary observation that 
$$ {1 \over {(\mu - \alpha)!}} = {{\mu^\alpha + p_{m-1}(\mu)} \over \mu !}  $$
where $|\alpha| = m$ and $p_{m-1}$ is of degree $m-1$, 
and on  resulting inequalities for large $ |\mu|$. See Corollary 1.

The second author acknowledges both the Institute for Advanced Study
 and the Mathematical Sciences Research Institute for support
during the academic year 1993-4. During that time he had productive
conversations with many mathematicians, but he particularly thanks Andy Nicas
for discussions about the mathematics in this paper. 
The first author acknowledges support from the NSF.

\noindent {\bf I. Preliminaries}

{\bf I.1 Notation and linear algebra}

We begin by introducing notation and recalling some elementary linear algebra.

We often use multi-index notation in this paper, usually without much comment.
If $\alpha$ is an $n$-tuple of non-negative integers, then we write $|\alpha| = \sum \alpha_j$ as usual.
 We write $\alpha ! = \prod \alpha_j !$.
As usual $z^\alpha = \prod z_j^{\alpha_j}$. Thus superscripts are powers, rather than being upper indices. 
We let $T$ denote the $n$-dimensional torus, and we let $d\theta$ denote the
invariant measure normalized so that $\int_T d\theta = 1$.
 We let $\theta = (\theta_1,...,\theta_n)$ be an element of $T$. 
Thus $ e^{i \alpha \theta}$ 
means $\prod e^{i \alpha_j \theta_j}$.

It will be important in this paper to
distinguish holomorphic polynomials from arbitrary polynomials. To help do so we
 sometimes write $h(z)$ when $h$ is holomorphic, and
 $f(z,{\overline z})$
when $f$ is arbitrary. 
Suppose that $g$ is a vector-valued holomorphic mapping. We write 
${\bf V}(g) = \{0\}$ to denote that the set of common
zeroes of the components of $g$ consists of the origin alone. The most important
 holomorphic mapping used here is $H_d$; as mentioned in
the introduction it is the tensor product of the identity mapping with itself $d$ times. In coordinates
 the components of $H_d$ are the
monomials $c_\alpha z^\alpha$ where $|\alpha|=d$ and $|c_\alpha|^2 = {d! \over {\alpha!}}$ of degree $d$.
Note that ${\bf V}(H_d) = \{0\}$. The components of $H_d$ are a basis for the complex vector
 space $V_d$ of (holomorphic) homogeneous polynomials of degree $d$ in $n$ variables, although we generally use
the monomials without these constants as the standard basis.

Suppose that $V$ is a finite-dimensional complex vector space, that $\zeta \in V$, and that
 $Q(\zeta,{\overline \zeta})=\sum_{\mu.\nu} Q_{\mu \nu} \zeta_\mu {\overline \zeta_\nu} $ is an Hermitian
form. By elementary linear algebra we can always find vectors $A_\mu$ and $B_{\mu}$ so that
$$  Q_{\mu \nu} = \langle A_\mu, A_\nu \rangle - \langle B_\mu, B_\nu \rangle. \eqno (1)$$ 
As a consequence of (1) we see that 

{\bf Lemma 1}. The Hermitian matrix $Q$ is positive definite if and only if
we have $ Q_{\mu \nu} = \langle A_\mu, A_\nu \rangle $ where the set of vectors $\{A_\mu\}$ form a basis. 

Lemma 1 implies the simple but fundamental
statement given by Lemma 2.

{\bf Lemma 2}. Let $f = \sum \sum c_{\alpha \beta} z^\alpha {\overline z}^\beta$ be a
 real-valued polynomial that is homogeneous of degree $m$ in both $z$ and ${\overline z}$.
Suppose that the underlying Hermitian matrix of coefficients $c_{\alpha \beta}$
has $k^+$ positive and $k^-$ negative eigenvalues.
 Then there are holomorphic vector-valued polynomials $P$ and $N$, each homogeneous of degree $m$, such that
$$  f(z,{\overline z}) = ||P(z)||^2 - ||N(z)||^2 . \eqno (2)$$
We can always choose $P$ and $N$ to have $k^+$ and $k^-$ components respectively.
The underlying matrix of coefficients is positive definite if and only if we may write
$$ f(z,{\overline z}) = ||P(z)||^2  $$
where  $ P = T H_{m+d}$ for some invertible linear transformation $T$. When the underlying matrix of coefficients
is positive definite we also must have   $ {\bf V}(P) = \{0\}$.

Proof. By elementary linear algebra as in (1) we may find vectors so that 
 $c_{\alpha \beta} = \langle A_\alpha, A_\beta \rangle - \langle B_\alpha, B_\beta \rangle$. Then we define $P$ by 
$P(z) = \sum A_\alpha z^\alpha $ and $N$ by $N(z) = \sum B_\alpha z^\alpha$. 
Easy computation then yields (2) and the result about the number of eigenvalues.
The eigenvalues are all positive if and only if $(c_{\alpha \beta})$ is positive definite.
Then the vectors $A_{\alpha}$ can be chosen to be linearly independent.  
Thus the underlying matrix is positive definite 
if and only if the matrix $(A_{\beta \alpha})$ is invertible,
where we write $A_\alpha = \sum A_{\beta \alpha} e_\beta $ in terms of the standard orthonormal basis. This implies
the statement about $T$ as well, because $H_{m+d}$ includes all monomials of degree $m+d$.
It remains to prove that ${\bf V}(P) = \{0\}$. 
Suppose that there exists $w$ with
$P(w)=0$. Then 
$$ 0 = ||P(w)||^2 = \sum \langle A_\alpha , A_\beta \rangle w^\alpha {\overline w}^\beta = 0 .$$
The underlying Hermitian form then has a non-trivial kernel  
 unless $w^\alpha$ must vanish for all $\alpha$. Hence, if it is positive definite, the variety must consist of the origin alone.
$\spadesuit$.

{\bf Example 1}. When $f(z,{\overline z}) = ||P(z)||^2$, and ${\bf V}(P) = \{0\}$, the underlying form 
might be only semi-definite. 
Consider the
simple example $ |z_1|^4 + |z_2|^4 $. The underlying form is diagonal, with eigenvalues $1,0,1$.
Theorem 1 implies however that we can obtain a positive definite underlying form  by multiplying by the squared norm. 

{\bf I.2 Positivity conditions}

Next we study positivity (away from the origin) of $f$ as a function. 
Suppose that $f$ is our given bihomogeneous polynomial, and that $f(z,{\overline z}) > 0 $ for $z \ne 0$.
Simple examples (see example III.2) reveal that the underlying matrix of coefficients can have eigenvalues of both
signs, so the Hermitian form on $V_m$ defined by
this matrix need not be positive. Nevertheless the
positivity of $f$ as a function has many consequences.

By compactness 
there is an $\epsilon > 0 $ such that $f(z,{\overline z}) \ge \epsilon$ on the unit sphere. By homogeneity we have
$$ f(z,{\overline z}) \ge \epsilon ||z||^{2m} \eqno (3)$$
everywhere.  
Now let the torus $T$ act on the sphere. This implies that, for all $\theta$, 
$$ g(\theta) = \sum c_{\alpha \beta} z^\alpha {\overline z}^\beta e^{i\theta (\alpha - \beta)} \ge \epsilon ||z||^{2m} \eqno (4)$$
as well. Integrating over the torus reveals that $\sum c_{\alpha \alpha} |z|^{2\alpha} \ge \epsilon ||z||^{2m}$ holds. 
To obtain more general inequalities we
 define $A_{\mu \nu}$ by
$$  A_{\mu \nu}(z,{\overline z}) = A_{\mu \nu} = 
\sum _{\beta} c_{(\beta + \mu - \nu)\beta} z^{\beta + \mu - \nu} {\overline z}^\beta $$
Note that $A_{\mu \nu}(z, {\overline z}) = \int_T e^{-i(\mu -\nu)\theta} g(\theta) d\theta$.

{\bf Proposition 1}. Suppose that (3) holds.
 For all non-zero $z$ the matrix $A_{\mu \nu}$ is positive definite and its minimum eigenvalue 
is at least $\epsilon ||z||^{2m}$.

Proof. We estimate  $ \sum A_{\mu \nu} \zeta_\mu {\overline \zeta_\nu}$ as follows.
$$ \sum A_{\mu \nu} \zeta_\mu {\overline \zeta_\nu} =
 \sum \int_T \zeta_\mu {\overline \zeta_\nu} e^{-i(\nu -\mu)\theta} g(\theta) d\theta $$
$$ = \int_T |\sum_{\mu} e^{-i\mu \theta} \zeta_\mu|^2 g(\theta) d\theta $$
$$ \ge  \int_T \epsilon ||z||^{2m} |\sum_{\mu} e^{-i\mu \theta} \zeta_\mu|^2  d\theta $$
$$ = \epsilon ||z||^{2m} \int_T |\sum_{\mu} e^{-i\mu \theta} \zeta_\mu|^2  d\theta $$
$$ = \epsilon ||z||^{2m} ||\zeta||^2 .$$
The inequality comparing the first and last expression proves the proposition. $\spadesuit$.

Next we make a simple observation that will enable us to relate these matrices to values of $f$.

{\bf Observation}. Let $\mu$ and $\alpha$ be multi-indices, and suppose that $|\alpha| = m$.
 Then there is a polynomial $p_{m-1}$ on ${\bf R^n}$ of degree $m-1$
 such that $ {1 \over {(\mu -\alpha)!}} - {\mu^{\alpha} \over {\mu !} } = {p_{m-1}(\mu) \over {\mu !}} $.
Hence there is a constant $C>0$ so that 
$ |{1 \over {(\mu -\alpha)!}} - {\mu^{\alpha} \over {\mu !} }| < C {|\mu|^{m-1} \over {\mu !}}$.

{\bf Proposition 2}. Given $\epsilon > 0$, there is an $N_0$ so that $|\mu| > N_0$ implies that
$$ |{1 \over {(\mu -\alpha)!}} - {\mu^{\alpha} \over {\mu !} }| < \epsilon |\mu|^m $$

Proof. This follows immediately from the observation and the definition of the limit, applied to 
$ {\rm lim}_{|\mu| \rightarrow \infty} {p_{m-1}(\mu) \over |\mu|^m} = 0 $. $\spadesuit$

{\bf I.3 Elementary Hilbert space considerations}

We now recall some facts about the Hilbert space of holomorphic square integrable functions
on the unit ball $B_n$. First of all, the Bergman kernel function $K(z,{\overline \zeta})$ for a bounded domain $\Omega$
is the
 kernel of the integral operator that projects the square integrable functions onto the closed subspace $A^2(\Omega)$ 
of holomorphic square integrable functions. It satisfies
 $ K(z,{\overline \zeta}) = \sum \phi_\alpha(z) {\overline \phi_\alpha(\zeta)} $
where $\{\phi_\alpha \}$ is a complete orthonormal set for $A^2(\Omega)$. For the unit ball
 (and a few other domains) normalized monomials form such a complete orthonormal
set and the series can be summed explicitly. 
For the unit ball $B_n$
the result is that $K(z,{\overline \zeta}) = {n! \over \pi^n} ( 1 - \langle z, \zeta \rangle)^{-(n+1)} $.
We can decompose $A^2 = A^2(B_n)$ as an orthogonal sum $A^2(B_n) = \sum V_k$ where $V_k$ is the subspace of 
holomorphic homogeneous polynomials of degree $k$. 

We say that a linear operator $L:L^2 \to L^2$ is positive if there is a positive number $c$ so that
$\langle Lg,g \rangle \ge c \langle g,g \rangle $. Here the inner product is the usual one given by integration.

{\bf Lemma 3}.  Suppose that $(c_{\alpha \beta})$ is the matrix of an Hermitian form on $V_m$.
Let $f(z,{\overline \zeta}) = \sum c_{\alpha \beta} z^\alpha {\overline \zeta}^\beta $ be the corresponding bihomogeneous
polynomial.
Let  $R: A^2 \to A^2$ be the operator defined by $$ Rg(z) = \int_{B_n} f(z,{\overline \zeta}) g(\zeta) dV(\zeta) .$$
Then $R$ is the zero operator on $V_m^{\perp}$ and $R$ maps $V_m$ to itself.  Furthermore  
$(c_{\alpha \beta})$ is positive definite if and only if the operator $R$ is positive on $V_m$.  

Proof. We expand $g$ as a convergent power series, writing $g(\zeta) = \sum g_\gamma \zeta^\gamma$.  We compute
$$ Rg(z) = \int \sum_{\alpha,\beta} \sum_{\delta} c_{\alpha \beta}  z^\alpha {\overline \zeta}^\beta g_\delta \zeta^\delta dV(\zeta).$$ 
By the orthogonality of the monomials in $A^2(B_n)$ we get zero unless $\beta = \delta$. This
yields 
$$ Rg(z) = \sum_{\alpha,\beta} c_{\alpha \beta} z^{\alpha} g_\beta ||z^{\beta}||_{A^2}^2  $$ 
and hence proves the first statement. Writing $d_\alpha = ||z^\alpha||^2_{A^2}$, we see that

$$ \langle Rg,g \rangle = \int \int \sum_{\alpha,\beta} \sum_{\gamma,\delta} c_{\alpha \beta} z^\alpha {\overline \zeta}^\beta
g_\delta \zeta^\delta {\overline g_\gamma} {\overline z}^\gamma dV(z) dV(\zeta) = 
 \sum_{\alpha, \beta} c_{\alpha \beta} d_\alpha d_\beta {\overline g_\alpha} g_\beta.$$

Thus, as a form on $V_m$, $R$ has underlying matrix $ (c_{ \beta \alpha} d_\alpha d_\beta) $. This is positive if and only if
$(c_{\alpha \beta})$ is, because of Lemma 1, and because the $d_\alpha$ are positive numbers. $\spadesuit$

We continue to let $f = \sum c_{\alpha \beta} z^\alpha {\overline \zeta}^\beta $ be our given bihomogeneous polynomial.
Let $R_{k+m}$ denote the operator on $A^2(B_n)$ whose kernel is $\langle z, \zeta \rangle^k f(z,{\overline \zeta})$.
The main point of Theorem 1 is to show that $f$ positive as a function
implies for sufficiently large $d$ that $||z||^{2d} f(z,{\overline z})$ 
has a positive definite underlying form.
By Lemma 3, it suffices to show that $R_{d+m}$ is a positive 
operator on $V_{m+d}$ for sufficiently large $d$.

\noindent {\bf II. Statement and proof of Theorem 1}

We are now prepared to give equivalent conditions for the positivity of a homogeneous polynomial on the sphere. 
Theorem 1 is the principal result of this paper.

{\bf Theorem 1}. Let $f(z,{\overline z}) = \sum \sum c_{\alpha \beta} z^\alpha {\overline z}^\beta$ be a
 real-valued polynomial
that is homogeneous of degree $m$ in $z$ and also in ${\overline z}$. The following are equivalent.

1) $f$ achieves a positive minimum value on the sphere.

2) There is an integer $d$ such that the underlying Hermitian matrix for $||z||^{2d} f(z,{\overline z})$ is positive
definite. Thus 
$$||z||^{2d} f(z,{\overline z}) = \sum E_{\mu \nu} z^\mu {\overline z}^\nu $$
where $(E_{\mu \nu})$ is positive definite.

3) There is an integer $d$ such that the operator $R_{m+d}$ defined by the kernel 
$ k_d(z,\zeta) = \langle z, \zeta \rangle ^d f(z, {\overline \zeta})$ is a positive operator from $V_m \subset A^2(B_n)$ to
itself.

4) There is an integer $d$ and a holomorphic homogeneous vector-valued 
polynomial $g$ of degree $m+d$ such that ${\bf V}(g) = \{0\}$ and such that
$||z||^{2d} f(z,{\overline z}) = ||g(z)||^2 $.

5) Write $f(z,{\overline z}) = ||P(z)||^2 - ||N(z)||^2 $ for holomorphic homogeneous vector-valued 
polynomials $P$ and $N$ of degree $m$. Then there is an integer $d$ and a linear transformation  $L$ such that the following are true:

5.1) $I - L^*L $ is positive semi-definite. 

5.2 ) $ H_d \otimes N = L (H_d \otimes P) $

5.3) $ {\bf V}(\sqrt{I-L^*L} (H_d \otimes P)) = \{0\}$

Proof. The hard implication is that 1) implies 3), so we do this last. It follows immediately from Lemma 3
that 2) and 3) are equivalent. We next observe that 2) implies 4),
immediately from Lemma 1. To prove that 2) implies 5), we recall that
 $||g\otimes h||^2 = ||g||^2 ||h||^2 $ and also that $ ||z||^{2d} = ||H_d(z)||^2 $. We plug these into 
$ ||g(z)||^2 = ||z||^{2d} (||P(z)||^2 - ||N(z)||^2)$, which is part of the hypothesis. 
After moving terms we obtain
$$ ||H_d \otimes N||^2 + ||g||^2 = ||H_d \otimes P||^2 .$$
Thus $ (H_d \otimes N) \oplus g $ and $ H_d \otimes P $ are holomorphic mappings with the same squared norms.
 Hence, see [D1], after 
perhaps including
enough zero components, there is a constant unitary matrix $U$ such that $( H_d \otimes N) \oplus g = U( H_d \otimes P) $. 
Writing this matrix  in terms of
blocks, one of which is $L$, we obtain submatrices $L$ and $\sqrt{I-L^* L}$ such that 5.1), 5.2), and 5.3) hold. 
Notice that we obtain the formula
that $ g = \sqrt{I-L^*L} (H_d \otimes P) $. Thus 4) and 5) are equivalent, without passing through 1).

That 5) implies 1) is also easy. Given 5), define $g$ by $ g = \sqrt{I-L^*L} (H_d \otimes P) $. 
Then 5.2) and the definition of $f$ imply that
$$ ||z||^{2d} ||N(z)||^2 + ||g(z)||^2 = ||z||^{2d} ||P(z)||^2 .$$
We see immediately that, on the sphere, $ ||N(z)||^2 + ||g(z)||^2 = ||P(z)||^2 $. 
Recalling that $f(z,{\overline z}) = ||P(z)||^2 - ||N(z)||^2 $, we see that the minimum value of $f$
on the sphere equals the minimum value of $g$
on the sphere, which is positive because ${\bf V}(g) = \{0\}$ by 5.3).

Thus the theorem follows if we can show that 1) implies 3).
Let $S$ be the operator on $A^2$ whose kernel is
$K(z,\zeta) f(z,\zeta)$. 
We claim that there are positive constants $c_{j,n}$ such that
 $c_{j,n} R_{j}$ equals the restriction of $S$ to $V_{j}$. 
The orthogonality of the spaces $V_j$ shows that $R_j$ is the zero operator off $V_j$.
Therefore the power series expansion  

$$ (1- \langle z,\zeta \rangle)^{-n-1} = \sum {{j+n} \choose j} \langle z, \zeta \rangle ^j $$
implies the claim.
To prove that $R_{d+m}$ is a positive 
operator on $V_{m+d}$ for
sufficiently large $d$, it therefore suffices to prove that the operator $S$
 has only finitely many non-positive eigenvalues. 

It remains only to 
show that $\langle Sg,g \rangle \ge c ||g||^2_{A^2}$ for $g \in V_{j}$ for sufficiently large $j$.
We introduce the operator $T$ on $L^2$, whose kernel is $K(z,\zeta) f(z,{\overline z})$. Then we have

$$ \langle Tg,g \rangle = \int \int K(z,\zeta) f(z, {\overline z}) {\overline g(z)}g(\zeta) dV(\zeta) dV(z) =
\int  f(z, {\overline z}) |g(z)|^2 dV(z) $$
by the reproducing property of the Bergman kernel. If we now introduce 
a cut-off function $\phi$ that is positive at the origin, and invoke the positivity of $f$ as a function,
we obtain that $\langle (T+\phi) g , g \rangle \ge c ||g||^2_{L^2}$.
Here $\phi$ is smooth, non-negative, positive at the origin,
and has compact support in the ball.

Next we write the operator equation $S = (T + \phi) +(E- \phi)$ where $E = S-T$. 
We observe that the kernel of $E$ vanishes at $z=\zeta$. Since this compensates for the
singularity of $K(z,\zeta)$ at $z = \zeta$ on the boundary of the ball,
$E$ is a compact operator. Since $\phi$ is smooth and compactly supported we also
have that $E-\phi$ is a compact operator.
Since $E - \phi$ is a compact operator on $L^2$, we can approximate it in norm by an operator $L$ with finite-dimensional
range. See [Ru] (pages 97-107) for the necessary statements about compact operators.
We write $|||L|||$ for the operator norm of $L$.

Then $S = (T + \phi) + (E- \phi - L) + L $.  Write $L_j$ for the restriction of $L$ to $V_j$.
Because $L$ has finite-dimensional range,
$ |||L_j||| < c/3 $ for $j$ sufficiently large.
Also we may assume that $ ||| E - \phi - L||| < c/3 $.  Therefore, on $V_j$ for sufficiently large $j$, the operator $S$ is a small
perturbation of $T + \phi$. Therefore $S$ is also a positive operator on $V_j$ for sufficiently large $j$. Hence $S$ has
only finitely many non-positive eigenvalues, so $R_{d+m}$ is a positive operator on $V_{m+d}$
for $d$ sufficiently large. This proves
that 1) implies 3) and completes the proof. $\spadesuit$.

The proof of Theorem 1 doesn't give a value for $d$; the value of $d$ is useful because it is closely related to
the embedding dimensions of mappings discussed in the next section, and because there
 is interest in the analogous
exponent in the classical theorem of Polya [R].
We therefore relate the numbers $E_{\mu \nu}$ to the numbers  ${c_{\alpha \beta}}$ and
 to the matrices
$A_{\mu \nu}$ from Proposition 1.

{\bf Corollary 1}.  Let
$ f(z,{\overline z}) =  \sum c_{\alpha \beta} z^\alpha {\overline z}^\beta$ 
 be a bihomogeneous polynomial that is positive away from
the origin. There is an integer $N$ such that the matrix $L_{\mu \nu}$ defined 
for $|\mu|= |\nu| = N $ by
$$  L_{\mu \nu} = \sum_{\alpha - \beta = \mu - \nu} c_{\alpha \beta}
 ({\mu \over |\mu|})^{\alpha} ({\nu \over |\nu|})^{\beta}  $$
is positive definite. Similarly there is an integer $N'$ so that the matrix 
$$  A_{\mu \nu} ( \sqrt{\nu \over |\nu|}, \sqrt {\nu \over |\nu|})  
 = \sum_{\alpha - \beta = \mu - \nu} c_{\alpha \beta}
 ({\nu \over |\nu|})^{\alpha + \beta }  $$
is positive definite when $ |\mu| = |\nu| = N'$.

Proof. For each positive integer $d$ we write $E_{\mu \nu}(d)$ for the underlying matrix of
coefficients for the function
 $$||z||^{2d} f(z, {\overline z}).$$
Then the multinomial expansion and interchange of summation yields

$$ E_{\mu \nu}(d) = d! \sum_{\alpha - \beta = \mu - \nu} c_{\alpha \beta} 
{1 \over (\mu - \alpha)! (\nu - \beta)!}.$$
We put $|\mu| = |\nu| = N$ and rewrite   
using the observation to get

$$ E_{\mu \nu}(d) = {N! \over N^m + p_{m-1}(N)} 
\sum_{\alpha - \beta = \mu - \nu} c_{\alpha \beta} 
(\mu^\alpha + p_{m-1}(\mu)) (\nu^\beta + p_{m-1}(\nu)) {1 \over \sqrt{ \mu ! \nu !}}$$
Next we multiply and divide by $N^m$ to get

$$ E_{\mu \nu}(d) = {N! \over \sqrt{ \mu ! \nu! } (1 + O(1/N))} 
\sum_{\alpha - \beta = \mu - \nu} c_{\alpha \beta} 
({ \mu \over |\mu|} + O(1/N)) ^\alpha  ({ \nu \over |\nu|} + O(1/N)) ^\beta  $$

Theorem 1 implies
that this matrix is positive definite when $d$ is sufficiently large. By Lemma 1
we do not alter the positive definiteness if we divide by 
${N! \over \sqrt{ \mu ! \nu!}}$, obtaining

$$D_{\mu \nu}(d) = {1 \over (1 + O(1/N))} 
\sum_{\alpha - \beta = \mu - \nu} c_{\alpha \beta} 
({ \mu \over |\mu|} + O(1/N)) ^\alpha  ({ \nu \over |\nu|} + O(1/N)) ^\beta  \eqno
(5)$$

Now it is clear that, for sufficiently large $|N|$, the matrix (5)
is positive definite if and only if 

$$ \sum_{\alpha - \beta = \mu - \nu} c_{\alpha \beta} 
({ \mu \over |\mu|}) ^\alpha  ({ \nu \over |\nu|}) ^\beta  $$
is positive. This proves the first statement. 

To prove the second, replace
$\mu$ by $ \mu = \nu + \alpha - \beta$ 
in the term $({ \mu \over |\mu|} + O(1/|N|)) ^\alpha$
arising in (5). Since the multi-indices $\alpha$ and $\beta$ have
fixed length $m$, we see that $ { {\mu + \alpha - \beta} \over |\mu|} $ can be also
written as $ {\nu \over |\nu|} + O(1/N) $. Thus the matrix 
$ A_{\mu \nu} ( \sqrt{\nu \over |\nu|}, \sqrt {\nu \over |\nu|})$ 
is also positive definite for sufficiently large $N$. 
This proves the Corollary. $\spadesuit$.

Observe that Proposition 1 guarantees that $ A_{\mu \nu} (z,{\overline z}) $ is
 positive definite for each
$z$. In Corollary 1 however the entries are  
evaluated at points depending on the indices. Corollary 1 enables one to obtain
quantitative information on the integer $d$ without explicitly performing the
multiplication by $||z||^{2d}$, but we do not develop this reasoning further here.

\noindent {\bf III. Examples and applications}

We first observe that the smallest possible value of the integer $d$ from Theorem 1 can be arbitrarily large.

{\bf Example 2}. Consider the one parameter
family of polynomials $p_c$ 
on ${\bf C^2}$ with variables $(z,w)$ defined by 
$$ p_c(z,w,{\overline z},{\overline w}) = |z|^4 + |w|^4 - c |zw|^2 .$$
Here $c$ is a real number. Then $p_c$ is positive on the sphere if and only if $ c < 2$. The original underlying matrix is diagonal, with
eigenvalues of $1,-c,1$. After multiplying through by $(|z|^2 + |w|^2)^d$ we also obtain diagonal matrices $L_d$ of
size $d+3$ by $d+3$. If we let $d(c)$
denote the minimum $d$ for which $L_d$ is positive definite, 
then we compute that $d(c) > (3c-2)/(2-c)$ and hence that
${\rm lim}_{c\rightarrow 2}d(c) = \infty$. 

This example indicates that the computations from Theorem 1 are easier when the original underlying matrix is diagonal.
 In fact this special case
is equivalent to the theorem of Polya [HLP]. Polya began with a homogeneous polynomial in several real variables,
 assumed to be positive
on the positive part of the hyperplane defined by $|x| = \sum x_j = 1$. He proved that there is an integer $d$ so that 
$$ |x|^d p(x) = q(x) $$
where all the coefficients of $q$ are positive. The value of this theorem is that $p$ is then expressed as 
a quotient of polynomials with positive
coefficients. When the underlying matrix of coefficients is diagonal, Theorem 1 implies Polya's theorem, and conversely. 
The correspondence
 between them is given by writing $x_j = |z_j|^2$ and observing that a diagonal matrix is positive definite if and
 only if all the diagonal
elements are positive. Our proof of Theorem 1 yields a new proof of Polya's theorem. On the other hand the computations
in Corollary 1 yield a direct proof, which 
we include now. 

{\bf Proof of Polya's theorem}. 

Put $f(x) = \sum c_\alpha x^\alpha$. Then 
$$|x|^d f(x) = d! \sum_\mu (\sum_\alpha {c_\alpha \over {(\mu -\alpha)!}}) x^\mu . $$
By the observation we can write the coefficient of $x^\mu$ as
$$ {d! \over {\mu!}} (f(\mu) + p_{m-1}(\mu)) $$
Suppose that $f(\mu) \ge \epsilon |\mu|^m $. By Proposition 2 we can choose $|\mu|$
 sufficiently large that $ |p_{m-1}(\mu)| < {\epsilon
\over 2} |\mu|^m $. Then $f(\mu) + p_{m-1}(\mu) > {\epsilon \over 2} |\mu|^m$ 
and the result follows. $\spadesuit$.

{\bf Remark 1}. The theorem fails in the positive semi-definite case. In other words,
the following is true. There exist real-valued polynomials $f(z,{\overline z})$
 that are non-negative as
 functions, and for which the underlying Hermitian form for
$||z||^{2d} f(z,{\overline z})$ has a negative eigenvalue for every $d$. One example is
$|z|^4 - 2 |zw|^2 + |w|^4 = (|z|^2 - |w|^2)^2 $. 

Next we give a concrete corollary obtained by applying the theorem when $n=2$. Suppose that $C= (c_{ab})$ 
is a Hermitian symmetric matrix of size $m+1$. 
We define a new matrix $E(C)$ by the following operation: $ E(C) = C \oplus 0 + 0 \oplus C$. This means that we augment
$C$ by including a leading row and column of zeroes, augment $C$ by including a last row and column of zeroes, and add the
two matrices to get a Hermitian symmetric matrix of size $m+2$. Then Theorem 1 gives a necesssary and
 sufficient condition for
there to be an integer $d$ for which $E^d (C)$ is positive definite. 

{\bf Corollary 2}.  There is an integer $d$ for which $E^d (C)$ is positive definite if and only if 
$\sum_{a=0}^m \sum _{b=0}^m c_{ab} t^a {\overline t}^b > 0 $ for all non-zero complex numbers $t$.

Proof. Replace $t$ by $z/w$ in the given condition, and multiply through by $|w|^{2m}$. The condition then becomes

$$ \sum_{a=0}^m \sum _{b=0}^m c_{ab} z^a w^{m-a} {\overline z}^b {\overline w}^{m-b} > 0 .$$

Thus we have a bihomogeneous polynomial in two variables. On the other hand, the operation $E$ is easily seen to be
equivalent to multiplying the function defining the underlying matrix of coefficients by $|z|^2 + |w|^2 $, and 
recomputing the underlying matrix of coefficients of the product. Thus the result follows from Theorem 1. $\spadesuit$

It may be worth remarking that the special case of Corollary 2 when the matrix is diagonal is a version of the Polya
theorem for inhomogeneous polynomials in one variable.

{\bf IV. Applications to proper mappings between balls}

We give several applications to proper mappings between balls. We consider two questions. First, suppose 
that $p: {\bf C^n} \to {\bf C^k}$ is a (holomorphic)
polynomial mapping and that $||p(z)|| < 1$ for all $z$ in the closed unit ball. We will show
that $p$ consists of some of the components of a proper mapping from $B_n$ to some $B_N$. 
We prove in Theorem 2
that there is a positive integer $N-k$ and a polynomial mapping $g:{\bf C^n} \to {\bf C^{N-k}}$ such that
$p \oplus g$ is the required mapping. In Theorem 4 we suppose that $q$ is a given polynomial function with $q(0)=0$ and that
$|q(z)| < 1$ on the closed ball. We prove that there is a rational proper mapping, reduced to lowest terms, with denominator
$1+q$.

{\bf Theorem 2}. Suppose that $p$ is a homogeneous
 vector-valued polynomial on $ {\bf C^n}$ and that $||p(z)||^2 < 1$ on the unit sphere. Then 
there is a polynomial mapping $g$
such that $p \oplus g$ defines a proper holomorphic mapping between balls.

We derive Theorem 2 from the following related result.

{\bf Theorem 3}. Suppose that $f(z,{\overline z})$ is a real-valued polynomial of degree $m$ (not necessarily homogeneous),
and that $f(z,{\overline z}) > 0$ on the unit sphere. Then there is a (holomorphic) polynomial mapping $g$ with finite-dimensional
range such that  $f(z,{\overline z}) = ||g(z)||^2 $ on the unit sphere.

Proof of Theorem 2 assuming Theorem 3. The result is obvious if $p$ is a constant.
Given a non-constant $p$, we consider $ f(z,{\overline z}) = ||z||^{2m} - p(z,{\overline z})$. Then $f$
satisfies the hypotheses of Theorem 3, so that $f(z,{\overline z}) = ||g(z)||^2 $ on the sphere for some $g$. Then
$p\oplus g$ is a non-constant holomorphic polynomial mapping whose squared norm $||p||^2 + ||g||^2$ equals unity on the sphere.
By the maximum principle $p\oplus g$ is the required mapping.

Proof of Theorem 3.

We are given that $f$ is degree $m$ in both $z$ and ${\overline z}$, so of total degree $2m$. We may assume that $m$ is even, 
because we may replace $f$ by $||z||^2 f$ otherwise without changing the hypothesis on $f$.
Given $f=  \sum_{|\alpha|} \sum_{|\beta|} c_{\alpha \beta} z^\alpha {\overline z}^\beta$,
we homogenize it by adding a new complex variable $t$, and put
$$ F(z,t,{\overline z},{\overline t}) =  
\sum_{\alpha} \sum_{\beta} c_{\alpha \beta} z^\alpha t^{m-|\alpha|} {\overline z}^\beta {\overline t}^{m-|\beta|}.$$
Since  $ F(z,e^{i \theta},{\overline z},e^{-i \theta}) = f(z e^{-i \theta},{\overline z e^{i \theta}})$, 
we see that
$F$ is positive when both $||z||^2 = 1$ and $|t|^2 = 1$. By homogeneity we then have

 $$ F(z,t,{\overline z},{\overline t}) \ge c (||z||^2  + |t|^2)^m $$
on the set $||z||^2 = |t|^2$. We can therefore choose a constant $C$ so that, away from the origin $(0,0)$, 

$$  F(z,t,{\overline z},{\overline t}) + C ((||z||^2  - |t|^2)^m  > 0 .$$

By Theorem 1, applied when the domain dimension is one larger, 
there is an integer $d$ and a holomorphic polynomial mapping $g$ such that 

$$ (||z||^2 + |t|^2)^d  F(z,t,{\overline z},{\overline t}) + C ((||z||^2  - |t|^2)^m = ||g(z,t)||^2 . \eqno (8) $$

Setting $t=1$ to dehomogenize (8), we see that 
$$ (||z||^2 + 1)^d  f(z,{\overline z}) + C ((||z||^2  - 1)^m = ||g(z,1)||^2 . \eqno (9)$$
It follows immediately from (9) that $2^d  f(z,{\overline z}) = ||g(z,1)||^2 $, so Theorem 3 follows. $\spadesuit$.

{\bf Remark}. Theorem 2 is easier when $p$ is homogeneous. One then applies Theorem 1 directly to 
the form by $f(z,{\overline z}) = ||z||^{2m} - ||p(z)||^2 .$
By statement 2) there is an integer $d$ and a
 holomorphic mapping $g$ such that
$$ ||g||^2 = ||z||^{2d} (||z||^{2m} - ||p(z)||^2) $$
The holomorphic mapping $g \oplus p$ then satisfies 
$$ ||p \oplus g||^2 = ||p||^2 + ||g||^2 = ||p||^2 + ||z||^{2d} ||p||^2 = ||z||^{2(m+d)} = ||H_{m+d}||^2 = 1 $$
on the unit sphere.

{\bf Remark}. Lempert [L1, L2] has proved an analog of Theorem 3 for functions extending holomorphically past the boundary of any
strongly pseudoconvex domain whose boundary is real-analytic. The difference here is that we
require that the mapping $g$ be both a polynomial
and be finite-dimensionally valued.

We consider briefly the consequences that $p \oplus g$ is a proper mapping between balls when $p$ is a 
vector-valued polynomial of degree $m$.
 Then $||p(z)||^2 + ||g(z)||^2 = 1$ on the sphere.
We write $p = \sum_{k=0}^m p_k$ for the decomposition of $p$ into its
homogeneous parts and similarly for $g$.
 The tensor product operation enables us to 
homogenize.  We have the equation (on the sphere) 
$$ \sum_{k,l=0}^{m+d} \langle g_k, g_l \rangle + \langle p_k, p_l \rangle = ||z||^{2(m+d)} .$$ 
Replace $z$ by $e^{i\theta} z$, where now $\theta$ is on the circle. We obtain the
following equations on the sphere:
$$ \sum (||g_k||^2 + ||p_k||^2) = ||z||^{2(m+d)} $$
when we consider the constant term, and 
$$ \sum ( \langle g_{l+s}, g_l \rangle + \langle p_{l+s}, p_{l} \rangle) = 0$$
for $s\not = 0$, when we equate Fourier coefficients. As in [D2] we can homogenize these equations by
tensoring each $g_k$ or $p_k$ with $H_{m+d-k}$. 
We obtain 
$$ \sum_k (||H_{m+d-k} \otimes g_k||^2 + ||H_{m+d-k} \otimes p_k||^2) = ||z||^{2(m+d)} =||H_{m+d}||^2 $$
from the constant term, and for $s \not = 0$, the equations (10).

$$ \sum_l ( \langle H_{m+d-l} \otimes  g_{l+s}, H_{m+d-l} \otimes g_l \rangle + 
\langle H_{m+d-l} \otimes q_{l+s}, H_{m+d-l} \otimes q_{l} \rangle)= 0 \eqno (10)$$

Hence the conclusion of Theorem 2 implies additionally that $g$ satisfies all these orthogonality relations.
In Proposition 1 from [D2], we are given an orthogonal direct sum of monomials of degree at most $m-1$.
We can find a monomial of degree $m$ that makes the entire expression into a proper mapping between balls, 
provided a certain form is positive semi-definite. Theorem 2 reveals that we do not need to assume
that the given monomial mappings are orthogonal in order to make them components of some proper polynomial
mapping. 

We next prove the theorem about allowable denominators.

{\bf Theorem 4}. Suppose that $q$ is a holomorphic polynomial that satisfies $q(0)=0$ and $|q|<1$
on the unit ball in ${\bf C^n}$. Then there is some $N$ and a polynomial mapping 
$p:{\bf C^n} \rightarrow {\bf C^N}$ such that ${p \over {1+q}}$ is a proper rational
mapping between balls that is reduced to lowest terms.

Proof.
Since $|q|<1$ on the sphere, for sufficiently small positive 
$\epsilon$, $|(1+\epsilon)q| < 1$ there as well. By Theorem 2 $(1+\epsilon)q$ 
is a component of a proper polynomial mapping $q \oplus g $ between balls. We consider an automorphism of
the target ball (See [D1] for example) of the form 
$\phi_a (w) = {{w - L_a(w)} \over {1 - \langle w, a \rangle}}$, where
$ a = {-1 \over {1 + \epsilon}} \oplus 0$. Then the composition $\phi_a (q \oplus g)$ has denominator $1+q$.
By considering the numerator it is easy to prove that the map constructed in this manner
is reduced to lowest terms, in the sense that $1+q$ does not divide all the components
of the numerator. $\spadesuit$.

Observe that the dimension $N$ depends on how close $q$ gets to
unity.

In this paper we have considered forms defined by polynomials. This enables us also 
to consider rational mappings.
The reader may wonder whether this may be too restrictive. Our results about proper mappings 
apply rather generally, because of a theorem of Forstneric. [F]. Suppose
that the domain dimension is at least two, and that $f$ is a proper holomorphic mapping between balls. If $f$ is smooth (of class
$C^\infty$) on the closed ball, then $f$ must be a rational mapping. 
For rational mappings ${p \over q}$ we consider the Hermitian
form defined by $||p||^2 - |q|^2$ and the methods of this paper apply. 
We hope to find further applications of these ideas to rational mappings.

Finally we remark on a geometric way to view Theorem 1. Consider the polynomial
mapping from ${\bf C^n}$ to ${\bf C^N}$ defined by
$\zeta_\alpha (z) = z^\alpha$. Here we assume that $N$ equals the dimension of the space of homogeneous polynomials of
degree $m$ in $n$ variables. The image of this mapping is an algebraic variety called the Veronese variety
$M_{m,n}$.
Our theorem has the following reformulation. Suppose that $Q$ is an Hermitian form on ${\bf C^N}$, and that $Q(\zeta,
{\overline \zeta})$ is positive for all $\zeta \in M_{m,n}$ except the origin. By example 2, 
$Q$ need not be positive everywhere.
There is however a Veronese variety $M_{d+m,n}$ of higher degree and an Hermitian form $E^d Q$ on it,
defined by the process of
multiplying the corresponding function by $||z||^{2d}$ and taking the underlying matrix of coefficients, such 
that the form is positive not only on $M_{d+m,n}$, but on the entire space (except the origin). A moment's thought
shows that both conditions are open. Thus the theorem remains true for varieties close to the Veronese. We can
think of the Veronese manifold as a testing manifold for positivity. This suggests that there is a generalization of
Theorem 1 to more general testing manifolds for positivity.

{\bf References}

[D1] D'Angelo, John P., Several Complex Variables and the Geometry of Real Hypersurfaces, CRC Press, Boca Raton, 1992.

[D2] D'Angelo, John P., The geometry of proper holomorphic maps between balls, Contemporary Mathematics, Volume 137, 1992, 191-215.

[F] Forstneric, F., Extending proper holomorphic mappings of positive codimension, Inventiones Math 95 (1989), 31-62.

[GH] Griffith, P., and Harris, J., Principles of Algebraic Geometry, John Wiley and Sons, New York, 1978.

[HLP] Hardy, G.H., Littlewood, J.E., and Polya,G., Inequalities, Cambridge, At the University Press, 1934.

[H] Handelman, David E., Positive Polynomials, Convex Integral Polytopes, and a Random Walk Problem, Lecture Notes in
Mathematics 1282, Springer-Verlag, Berlin.

[L1] Lempert, L., Imbedding Cauchy-Riemann manifolds into a sphere, International Journal of Math 1 (1990), 91-108.

[L2] Lempert, L., Imbedding pseudoconvex domains into a ball, American Journal of Math 104 (1982), 901-904.

[R] Reznick, Bruce, Universal denominators in Hilbert's seventeeth problem, Math Z. (1995), (to appear).

[Ru] Rudin, Walter, Functional Analysis, McGraw-Hill, New York, 1973.

\end